\documentclass[12pt]{article}
\usepackage{amssymb,amsmath}

\newcommand{\Z}{{\mathbb Z}}
\newcommand{\Q}{{\mathbb Q}}
\newcommand{\Fp}{{\mathbb F}_p}
\newcommand{\Fps}{{\mathbb F}_{p^2}}
\newcommand{\Fpn}{{\mathbb F}_{p^{\nu}}}
\newcommand{\Zpn}{{\mathbb Z}_{p^{\nu}}}
\newcommand{\Zps}{{\mathbb Z}_{p^2}}

\newcommand{\OO}{\mathcal{O}}
\newcommand{\M}{\mathcal{M}}
\renewcommand{\char}{\text{char}}
\newcommand{\Gal}{\text{Gal}}
\newcommand{\Res}{\text{Res}}

\newcommand{\N}{{\rm N}}
\newcommand{\Tr}{\text{Tr}}
\newcommand{\Kbar}{\overline{K}}

\newcommand{\cbar}{\overline{c}}
\newcommand{\proof}{\noindent{\em Proof: }}
\newcommand{\qed}{\hspace{\fill}$\square$}
\newcommand{\ra}{\rightarrow}
\newcommand{\lra}{\longrightarrow}
\newcommand{\dst}{\displaystyle}
\newcommand{\mubold}{{\boldsymbol\mu}}

\newtheorem{theorem}{Theorem}
\newtheorem{lemma}[theorem]{Lemma}
\newtheorem{prop}[theorem]{Proposition}
\newtheorem{cor}[theorem]{Corollary}
\newenvironment{algorithm}{\noindent\refstepcounter{theorem}{\bf Algorithm \thesection.\arabic{theorem}} }{}
\numberwithin{equation}{section}
\numberwithin{theorem}{section}

\title{Indices of inseparability for elementary abelian
$p$-extensions}
\author{Kevin Keating \\
Department of Mathematics \\
University of Florida \\
Gainesville, FL 32611 \\
USA \\[.2cm]
{\tt keating@ufl.edu}}

\begin{document}

\maketitle

\begin{abstract}
\noindent
Let $K$ be a local field whose residue field $\Kbar$ is
a finite field of characteristic $p$ and let $L/K$ be a
finite totally ramified Galois extension.  Fried
\cite{fried} and Heiermann \cite{heier} defined the
``indices of inseparability'' of $L/K$, a refinement of
the ramification data of $L/K$.  We give a method for
computing the indices of inseparability of the extension
$L/K$ in terms of the norm group $\N_{L/K}(L^{\times})$
in the case where $K$ has characteristic $p$ and
$\Gal(L/K)$ is an elementary abelian $p$-group with a
single ramification break.  In some cases our methods
lead to simple formulas for the indices of inseparability.
\end{abstract}

\section{Introduction}

     Let $K$ be a local field with finite residue field
$\Kbar$ and let $L/K$ be a finite separable extension.
The {\em indices of inseparability} of $L/K$ were defined
by Fried \cite{fried} for fields of characteristic $p$,
and by Heiermann \cite{heier} for fields of characteristic
0.  The indices of inseparability of $L/K$ constitute a
refinement of the usual ramification data for $L/K$, as
described for instance in Serre \cite[IV]{cl}, in that
the indices of inseparability determine the ramification
data, but the ramification data does not always
determine the indices of inseparability.  More
precisely, suppose $L/K$ is Galois, with Galois group
$G=\Gal(L/K)$.  Then $G$ has a filtration by the lower
ramification subgroups $G\ge G_0\ge G_1\ge\dots$, as
defined in \cite[IV\,\S1]{cl}, and the ramification data
of $L/K$ is determined by the sequence $(|G_i|)_{i\ge0}$.
Each ramification subgroup $G_i$ is normal in $G$, and
the quotients $G_i/G_{i+1}$ for $i\ge1$ are elementary
abelian $p$-groups.  If these quotients are all cyclic
then the ramification data of $L/K$ determines the
indices of inseparability, but if there is some $i\ge1$
such that the quotient $G_i/G_{i+1}$ is not cyclic then
the indices of inseparability cannot be completely
determined from the orders of the groups $G_i$.

     Let $K^{alg}$ be an algebraic closure of $K$.  Class
field theory gives a one-to-one correspondence between
the set of finite abelian subextensions $L/K$ of
$K^{alg}/K$ and the set of closed finite-index subgroups
$H$ of $K^{\times}$.  This correspondence maps the
extension
$L/K$ onto the group $\N_{L/K}(L^{\times})$, where
$\N_{L/K}:L\ra K$ is the norm map.  In principle it should
be possible to determine the indices of inseparability
of $L/K$ in terms of the norm group
$\N_{L/K}(L^{\times})$.  In fact it is straightforward
to compute the ramification data of
$L/K$ in terms of $\N_{L/K}(L^{\times})$ (see for instance
Theorem 2 in \cite[XV\,\S2]{cl}), but extracting the
additional information contained in the indices of
inseparability seems to be considerably more difficult.
Perhaps this should not be surprising, since the
ramification data of $L/K$ can be defined in terms of
the norm $\N_{L/K}$.  This is the approach taken by
Fesenko and Vostokov in \cite[III\,\S3]{FV}, but their
definition does not extend
in any obvious way to give a general norm-based
definition of the indices of inseparability.

     In this paper we consider a special case of the
problem of determining the indices of inseparability of
an abelian extension in terms of its norm group.  In
order to isolate the difficulties involved in passing
from the usual ramification data to the refinement given
by the indices of inseparability we restrict our
attention to the case of an extension with a single
ramification break.  Let $L/K$ be a finite totally
ramified abelian extension whose Galois group
$G=\Gal(L/K)$ has a single ramification break $b\ge1$.
Then $G\cong G_b/G_{b+1}$ is an elementary abelian $p$-group.
Thus if $[L:K]>p$ we are in a situation where knowledge of
the ramification data does not determine the indices of
inseparability.

     In section \ref{ram} we define the indices of
inseparability of a finite totally ramified extension
of local fields $L/K$.  In section \ref{Ab} we give an
explicit computation of the norm group
$\N_{L/K}(L^{\times})$ in the case where $L/K$ is an
elementary abelian $p$-extension with a single
ramification
break.  In section \ref{AS} we assume $\char(K)=p$ and
interpret the results of section \ref{Ab} in terms of
Artin-Schreier theory.  This leads to an algorithm for
computing the indices of inseparability in terms of the
Artin-Schreier equations which define $L$.  In section
\ref{explicit} we use the results of sections \ref{Ab}
and \ref{AS} to give explicit formulas for the indices
of inseparability of certain elementary abelian
$p$-extensions in characteristic $p$.
In section \ref{char0} we show how the results of
section \ref{explicit} can be extended to
apply to some extensions of local fields of
characteristic 0.

     The author would like to thank the referee for
several suggestions which improved the exposition, and
for an unusually thorough and careful reading of the
paper.

\section{Ramification and indices of inseparability}
\label{ram}

In this section we define the ramification data and
indices of inseparability for a finite totally ramified
Galois extension of local fields.  We then describe
how the indices of inseparability are computed in
terms of the ramification data.

     Let $K$ be a local field with finite residue field
$\Kbar$ of characteristic $p$ and let $v_K$ be the
normalized valuation on $K$.  Let $e=v_K(p)$ be the
absolute ramification index of $K$; thus $e=\infty$ if
and only if $\char(K)=p$.  Let
$\OO_K=\{\alpha\in K:v_K(\alpha)\ge0\}$ be the ring of
integers of $K$, let $\M_K$ be the maximal ideal of $\OO_K$,
and for $r\ge1$ let $U_K^r=1+\M_K^r$.  Let $L/K$ be a finite
separable totally ramified Galois subextension of $K^{alg}/K$
of degree $n=ap^{\nu}$, with $p\nmid a$.  Also
let $v_p$ denote the $p$-adic valuation on $\Z$.

   We begin by recalling the definition of the
ramification data of $L/K$.  Set $G=\Gal(L/K)$ and let
$\pi_L$ be a uniformizer for $L$.  For a nonnegative
real number $x$ we define the $x$th ramification group
of $G$ (with respect to the lower numbering) by
\begin{equation}
G_x=\{\sigma\in G:v_L(\sigma(\pi_L)-\pi_L)\ge x+1\}.
\end{equation}
We easily deduce the following: $G_x$ is a normal
subgroup of $G$ which is independent of the choice of
$\pi_L$; $G_x\le G_y$ for $y\le x$; $G_0=G$; $G_x=\{1\}$
for all sufficiently large $x$; and $G_x=G_i$, where
$i=\lceil x\rceil$ is the smallest integer such that
$i\ge x$.

     Let $b\ge0$.  If there is $\sigma\in G$ such that
$v_L(\sigma(\pi_L)-\pi_L)=b+1$ then we say that $b$ is a
lower ramification break for $G$; in this case $b$ must
be an integer.  This is equivalent to having $b\in\Z$
and $G_{b+1}\not=G_b$.  If $b\ge1$ then there is a group
embedding of $G_b/G_{b+1}$ into $\M_L^b/\M_L^{b+1}$
which carries $\sigma G_{b+1}$ onto
$\dst\frac{\sigma(\pi_L)-\pi_L}{\pi_L}+\M_L^{b+1}$.
Hence $G_b/G_{b+1}$ is an elementary abelian $p$-group.

     We now recall the definition of the indices of
inseparability $i_0,i_1,\dots,i_{\nu}$ of $L/K$ as formulated
by Fried in the case $\char(K)=p$ \cite[pp.\,232--233]{fried}
(see also \cite[\S2]{fm}), and by Heiermann in the case
$\char(K)=0$ \cite[\S3]{heier}.  Set $q=|\Kbar|$
and let
$\mubold_{q-1}$ denote the group of $q-1$ roots of unity
in $K$.  Then $R=\mubold_{q-1}\cup\{0\}$ is the set of
Teichm\"uller representatives for $K$; if $\char(K)=p$ then
$R$ is a subfield of $K$ which can be identified with the
residue field $\Kbar$.  Given uniformizers $\pi_K$,
$\pi_L$ for $K$, $L$ there are unique $c_h\in R$ such that
$\dst\pi_K=\sum_{h=0}^{\infty}c_h\pi_L^{h+n}$.  Let
$0\le j\le\nu$ and set
\begin{equation}
\tilde{\imath}_j=
\min\{h\ge0:c_h\not=0,\,v_p(h+n)\le j\};
\end{equation}
if $c_h=0$ for all $h\ge0$ such that
$v_p(h+n)\le j$ let $\tilde{\imath}_j=\infty$.
Note that since $v_p(n)=\nu$ we have
$\tilde{\imath}_{\nu}=0$ and
\begin{equation}
\tilde{\imath}_j=\min\{h\ge0:c_h\not=0,\,v_p(h)\le j\}
\end{equation}
for $0\le j<\nu$.
The indices of inseparability of $L/K$ are now defined
recursively by $i_{\nu}=\tilde{\imath}_{\nu}=0$ and
$i_j=\min\{\tilde{\imath}_j,i_{j+1}+ne\}$ for
$j=\nu-1,\dots,1,0$.  Equivalently,
\begin{equation} \label{ij}
i_j=\min\{\tilde{\imath}_{j'}+ne(j'-j):j\le j'\le\nu\}.
\end{equation}
If $\char(K)=0$ then
$\tilde{\imath}_j$ may depend on the choice of $\pi_L$
(but not on the choice of $\pi_K$).  Nevertheless,
$i_j$ is independent of the choice of uniformizers
\cite[Th.\,7.1]{heier}.  If $\char(K)=p$ then it is
easily seen that $i_j=\tilde{\imath}_j$ does not depend
on the choice of uniformizers.  It follows immediately
from the definitions that
$0=i_{\nu}\le i_{\nu-1}\le\dots\le i_0$.
Furthermore, if $v_p(i_j)=j'<j$ then
$i_j=i_{j-1}=\dots=i_{j'}$.

     The connection between the ramification data of the
extension $L/K$ and the indices of inseparability of $L/K$
can be described most easily in terms of the Hasse-Herbrand
function of $L/K$, which is defined by
\begin{equation} \label{phi}
\phi_{L/K}(x)=\int_0^x\frac{dt}{|G_0:G_t|}
\end{equation}
for $x\ge0$.  Thus $\phi_{L/K}:[0,\infty)\ra[0,\infty)$
is a continuous increasing piecewise linear function
which is differentiable at all positive values except for
the lower ramification breaks of $L/K$.  If $x>0$ is not
a ramification
break then the derivative of $\phi_{L/K}$ at $x$ is
$\phi_{L/K}'(x)=|G_x|/n$, where $n=ap^{\nu}=[L:K]=|G_0|$.
Hence the ramification data of $L/K$ is encoded in the
function $\phi_{L/K}$.

     One can recover the Hasse-Herbrand function $\phi_{L/K}$
(and hence the orders of the ramification subgroups)
from the indices of inseparability of $L/K$ by the formula
\begin{equation} \label{recover}
\phi_{L/K}(x)=\frac1n\cdot\min\{i_j+p^jx:0\le j\le\nu\}
\end{equation}
(see \cite[Cor.\,6.11]{heier}).  Let $S$ be one of the
segments that make up the graph of $\phi_{L/K}$.  It
follows from (\ref{recover}) that $S$ can be extended to
meet the $y$-axis at $(0,i_j/n)$ for some $0\le
j\le\nu$.  It follows
that if $L/K$ has $k$ positive ramification breaks then
the graph of $\phi_{L/K}$ determines $k+1$ indices of
inseparability for $L/K$.  In particular, if $k=\nu$ then
the ramification data of $L/K$
determines all $\nu+1$ indices of inseparability.  The
difficult cases where $L/K$ has fewer than $\nu$ positive
ramification breaks occur when $L/K$ has a positive
ramification break $b$ such that the elementary abelian
$p$-group $G_b/G_{b+1}$ has rank greater than 1.

\section{Abelian extensions with one break} \label{Ab}

Let $K$ be a local field whose residue field $\Kbar$ is
a finite field of characteristic $p$, and let $L/K$ be a
nontrivial totally ramified Galois extension with a
single ramification break $b\ge1$.
We assume for simplicity that if $\char(K)=0$ then $p>2$
and $p\nmid b$.  (The condition $p\nmid b$ is automatic
if either $\char(K)=p$ or $[L:K]>p$.) The aim of this
section is to give an explicit computation of the norm
group $H=\N_{L/K}(L^{\times})$.  In section \ref{AS} we
will use this computation in the case $\char(K)=p$ to
get information about the subgroup
of the additive group of $K$ which corresponds to $L/K$
under Artin-Schreier theory.  We note that Monge uses
the same methods in Chapter 6 of \cite{monge} to compute
the norm groups of certain extensions of local fields of
characteristic 0.

     Set $G=\Gal(L/K)$.  Then $G_b=G$ and $G_{b+1}=\{1\}$,
so $G\cong G_b/G_{b+1}$ is an elementary abelian $p$-group,
say $G\cong(\Z/p\Z)^{\nu}$ with $\nu\ge1$.  It follows
from (\ref{phi}) that the graph of the function
$\phi_{L/K}(x)$ consists of a line segment of slope 1
from $(0,0)$ to $(b,b)$ and a ray of slope $p^{-\nu}$
starting at $(b,b)$.  Hence by (\ref{recover}) we have
$i_j\ge bp^{\nu}-bp^j$ for $0\le j\le\nu$, with
equality for $j=0$ and $j={\nu}$.

     Let $E_K(X)=\prod_{p\nmid c}\,(1-X^c)^{-\mu(c)/c}
\in\OO_K[[X]]$ denote the Artin-Hasse exponential series
(cf.~\cite[III\,\S1]{pdg}), where $\mu$ is the M\"obius
function.  Of course, $E_K(X)$ depends only on $p$ and
on the characteristic of $K$: If $\char(K)=p$ then
$E_K(X)=E_p(X)\in\Fp[[X]]$, while if $\char(K)=0$ then
$E_K(X)=E_0(X)\in\Z_p[[X]]$.
Since $E_K(X)=1+X+\dots$, the norm group of
$L/K$ can be described in terms of this series:

\begin{prop} \label{generate}
Let $\pi_L$ be a uniformizer for $L$.  The subgroup
$H=\N_{L/K}(L^{\times})$ of $K^{\times}$ is generated by
$U_K^{b+1}$, $(K^{\times})^p$, $\N_{L/K}(\pi_L)$, and
the set
\begin{equation}
\{\N_{L/K}(E_K(r\pi_L^k)):r\in R,\;
1\le k\le b,\;p\nmid k\}.
\end{equation}
\end{prop}

\proof Since $b$ is the only lower ramification break of
$L/K$, it follows from \cite[V\,\S6, Cor.\,3]{cl} that
$\N_{L/K}(U_L^{b+1})=U_K^{b+1}$.  Since
$K^{\times}/H\cong G$ is killed by $p$ we have
$H\supset(K^{\times})^p$, so $H$ contains all the
listed elements.  Since $p\nmid b$ the group
$L^{\times}$ is generated by $U_L^{b+1}$,
$(L^{\times})^p$, $\pi_L$, and the set
\begin{equation}
\{E_K(r\pi_L^k):r\in R,\;1\le k\le b,\;p\nmid
k\}.
\end{equation}
Therefore $H$ is generated by the listed elements.
\qed \medskip

     Let $\pi_L$ be a uniformizer of $L$ and let
\begin{equation}
g(X)=X^{p^{\nu}}+a_1X^{p^{\nu}-1}+\dots
+a_{p^{\nu}-1}X+a_{p^{\nu}}
\end{equation}
be the minimum polynomial for $\pi_L$ over $K$.  For
$0\le j\le\nu$ define
\begin{equation} \label{ak}
\widehat{\imath}_j=\min\{p^{\nu}v_K(a_k)-k:
1\le k\le p^{\nu},\;v_p(k)\le j\}.
\end{equation}
Then by \cite[Prop.\,3.12]{heier} we have
\begin{equation} \label{hatij}
i_j=\min\{\widehat{\imath}_{j'}+p^{\nu}e(j'-j):j\le j'\le\nu\}
\end{equation}
(cf.\ (\ref{ij})).  In particular, if $\char(K)=p$ then
$i_j=\widehat{\imath}_j$.
Set $t=\N_{L/K}(\pi_L)$; then $t$ is a uniformizer for
$K$, and it follows from our assumptions that
$t=(-1)^{p^{\nu}}a_{p^{\nu}}=-a_{p^{\nu}}$.  For $k\ge1$
write $k=k_0+k_1p^{\nu}$ with $1\le k_0\le p^{\nu}$ and
define $a_k=t^{k_1}a_{k_0}$.  Since $i_0=p^{\nu}b-b$ and
$p\nmid b$, by (\ref{hatij}) we have
$i_0=\widehat{\imath}_0$.  By (\ref{ak}) we get
$v_K(a_{b_0})=b-b_1$, and hence $v_K(a_b)=b$.

     Recall that $R=\mubold_{q-1}\cup\{0\}$ is the set
of Teichm\"uller representatives of $K$.  There are
unique $c_{i,h}\in R$ such that
$\dst a_i=\sum_{h=1}^{\infty}\,c_{i,h}t^h$.  It follows
from (\ref{ak}) that for $0\le j\le\nu$ we have
\begin{align}
\widehat{\imath}_j
=\min\{hp^{\nu}-i:1\le i\le p^{\nu},\;v_p(i)\le j,
\;h\ge1,\;c_{i,h}\not=0\}.
\end{align}
Since $a_{i+p^{\nu}}=ta_i$ we see that $c_{i,h}$ depends
only on the value of $hp^{\nu}-i$.  It follows that
\begin{align} \label{without}
\widehat{\imath}_j&=\min\{hp^{\nu}-i:i\ge1,\;v_p(i)\le j,\;h\ge1,\;c_{i,h}\not=0\}.
\end{align}
Using the fact that
$bp^{\nu}-bp^j\le i_j\le i_0=bp^{\nu}-b$ we get
\begin{align} \label{bpn}
\widehat{\imath}_j&=\min\{bp^{\nu}-i:b\le i\le bp^j,\;v_p(i)\le j,
\;c_{i,b}\not=0\}.
\end{align}

     For $k\ge1$ such that $p\nmid k$ set
\begin{equation} \label{gX}
g_k(X)=(-1)^{k+1}\prod_{\zeta^k=1}\,g(\zeta X^{1/k}),
\end{equation}
where the product is taken over the $k$th roots of unity
$\zeta\in K^{alg}$.  Then $g_k(X)$ is a monic polynomial
of degree $p^{\nu}$ with coefficients in $\OO_K$ whose
roots are the $k$th powers of the roots of $g(X)$.  For
$r\in R$ set $h_k^r(X)=-r^{p^{\nu}}\cdot g_k(-r^{-1}(X-1))$.
Then $h_k^r(X)$ is a polynomial of degree
$p^{\nu}$ with coefficients in $\OO_K$ whose multiset of
roots is $\{1-r\omega^k:g(\omega)=0\}$.  Since we are
assuming either $\char(K)=p$ or $p>2$ we see that
$h_k^r(X)$ is monic and
\begin{equation} \label{umr}
\N_{L/K}(1-r\pi_L^k)=-h_k^r(0)=r^{p^{\nu}}\cdot
g_k(r^{-1}).
\end{equation}

\begin{lemma} \label{fi}
For $i\ge1$ we have $v_K(a_i)\ge f_i$, where 
$\dst f_i=\left\lceil b-\frac{p^{v_p(i)}b-i}{p^{\nu}}
\right\rceil$. 
\end{lemma}

\proof Let $j=v_p(i)$.  If $j\ge\nu$ then
$a_i=-t^{i/p^{\nu}}$ and the claim is obvious.  If
$j<\nu$ then we may assume $1\le i\le p^{\nu}$.
It follows from (\ref{ak}) and the inequality
$\widehat{\imath}_j\ge i_j\ge bp^{\nu}-bp^j$
that
\begin{equation}
v_K(a_i)\ge p^{-\nu}(\widehat{\imath}_j+i)\ge
b-bp^{j-\nu}+p^{-\nu}i.
\end{equation}
Hence $v_K(a_i)\ge f_i$ in this case as well.
\qed \medskip

     Let
$\tilde{g}(X)=a_1X^{p^{\nu}-1}+\dots+a_{p^{\nu}-1}X$;
then $g(X)=X^{p^{\nu}}+\tilde{g}(X)-t$.
Choose $s\in K^{alg}$ such that $s^k=r$.  It follows
from (\ref{gX}) and (\ref{umr}) that
\begin{equation}
\N_{L/K}(1-r\pi_L^k)=(-1)^{k+1}r^{p^{\nu}}\cdot
\left(\prod_{\zeta^k=1}
(\zeta^{p^{\nu}}s^{-p^{\nu}}-t+\tilde{g}(\zeta
s^{-1}))\right).
\end{equation}
By Lemma~\ref{fi} we have
$v_K(a_i)\ge\lceil(1-p^{-1})b\rceil$ for
$1\le i\le p^{\nu}-1$.  Since
$2\lceil(1-p^{-1})b\rceil\ge b+1$ we get the following
congruences modulo $\M_K^{b+1}$:
\begin{alignat}{2}
\N_{L/K}(1-r\pi_L^k)&\equiv
(-1)^{k+1}r^{p^{\nu}}\cdot\left(\prod_{\zeta^k=1}
(\zeta^{p^{\nu}}s^{-p^{\nu}}-t)\right)\cdot
\left(1+\sum_{\zeta^k=1}
\frac{\tilde{g}(\zeta
s^{-1})}{\zeta^{p^{\nu}}s^{-p^{\nu}}-t}\right)
\\
&\equiv(1-r^{p^{\nu}}t^k)\cdot\left(1+\sum_{\zeta^k=1}
\frac{\zeta^{-p^{\nu}}s^{p^{\nu}}\tilde{g}(\zeta
s^{-1})}{1-ts^{p^{\nu}}\zeta^{-p^{\nu}}}
\right)\\
&\equiv(1-r^{p^{\nu}}t^k)\;\cdot 
\left(1+\sum_{\zeta^k=1}
\sum_{i=1}^{p^{\nu}-1}\sum_{j=0}^{\infty}\,
t^ja_i\zeta^{-jp^{\nu}-i}s^{jp^{\nu}+i}\right).
\end{alignat}
Since $\sum_{\zeta^k=1}\zeta^{-jp^{\nu}-i}$ equals $k$
if $k\mid jp^{\nu}+i$, and 0 otherwise, we get the
following congruences modulo $\M_K^{b+1}$:
\begin{alignat}{2}
\N_{L/K}(1-r\pi_L^k)&\equiv(1-r^{p^{\nu}}t^k)\cdot
\left(1+k\cdot\sum_{j=0}^{\infty}
\sum_{\substack{0<i<p^{\nu}\\[.05cm]k\mid jp^{\nu}+i}}
t^ja_is^{jp^{\nu}+i}\right)\\
&\equiv(1-r^{p^{\nu}}t^k)\cdot
\left(1+k\cdot\sum_{\substack{k\mid m\\p^{\nu}\nmid
m}}\,a_mr^{m/k}\right) \\
&\equiv(1-r^{p^{\nu}}t^k)\cdot
\left(1+k\cdot\sum_{p^{\nu}\nmid l}\,a_{lk}r^l\right).
\label{above}
\end{alignat}

     In order to apply Proposition~\ref{generate} we
express our norm computations in terms of the
Artin-Hasse exponential:
\begin{alignat}{2}
\N_{L/K}(E_K(r\pi_L^k))
&=\prod_{p\nmid c}\,\N_{L/K}(1-r^c\pi_L^{ck})^{-\mu(c)/c}
\\
&\equiv\prod_{p\nmid c}\,(1-r^{cp^{\nu}}t^{ck}
)^{-\mu(c)/c}\cdot 
\prod_{p\nmid c}\left(1+ck\cdot\sum_{p^{\nu}\nmid l}
\,a_{lck}r^{lc}\right)^{-\mu(c)/c}\\
&\equiv E_K(r^{p^{\nu}}t^k)\cdot
\left(1-k\cdot\sum_{p\nmid c}\sum_{p^{\nu}\nmid l}\,
\mu(c)a_{lck}r^{lc}\right),
\end{alignat}
with all congruences modulo $\M_K^{b+1}$.
By setting $m=lc$ we get
\begin{alignat}{2}
\N_{L/K}(E_K(r\pi_L^k))
&\equiv E_K(r^{p^{\nu}}t^k)\cdot
\left(1-k\cdot\sum_{p\nmid c}\sum_{\substack{c\mid
m\\p^{\nu}\nmid m}}\,
\mu(c)a_{mk}r^{m}\right)\\
&\equiv E_K(r^{p^{\nu}}t^k)\cdot
\left(1-k\cdot\sum_{p^{\nu}\nmid
m}a_{mk}r^{m}\sum_{\substack{c\mid
m\\p\nmid c}}\,
\mu(c)\right)\\
&\equiv E_K(r^{p^{\nu}}t^k)\cdot
\left(1-k\cdot\sum_{j=0}^{\nu-1}\,a_{kp^j}r^{p^j}\right)
\\
&\equiv E_K(r^{p^{\nu}}t^k)\cdot
E_K\left(-k\cdot\sum_{j=0}^{\nu-1}\,a_{kp^j}r^{p^j}\right).
\label{normE}
\end{alignat}

     For $1\le k\le b$ and $0\le j\le\nu-1$ define
\begin{equation}
S_{kj}=\{h\in\Z:f_{kp^j}\le h\le b\}.
\end{equation}
Using the expansion
$a_i=\dst\sum_{h=1}^{\infty}\,c_{i,h}t^h$ we can rewrite
(\ref{normE}) in the form
\begin{alignat}{2} \label{NEp}
\N_{L/K}(E_K(r\pi_L^k))&\equiv
E_K(r^{p^{\nu}}t^k)
\cdot E_K\left(-k\cdot\sum_{j=0}^{\nu-1}\sum_{h\in S_{kj}}
c_{kp^j,h}r^{p^j}t^h\right)\\
&\equiv E_K(r^{p^{\nu}}t^k)\cdot\prod_{j=0}^{\nu-1}
\prod_{h\in S_{kj}}E_K(-kc_{kp^j,h}r^{p^j}t^h).
\label{prod}
\end{alignat}
As before all congruences are taken modulo $\M_K^{b+1}$.

\section{Artin-Schreier theory} \label{AS}

Let $K$  be a local field of characteristic $p$ with
finite residue field $\Kbar$, and let $L/K$ be a finite
totally ramified Galois extension with a single
ramification break $b\ge1$.  As above $\pi_L$ is a
uniformizer for $L$ and $t=\N_{L/K}(\pi_L)$ is a
uniformizer for $K$.  In this section we interpret the
norm computations from section~\ref{Ab} in terms of
Artin-Schreier theory.

     Let $K^{ab}/K$ denote the maximum abelian subextension
of $K^{alg}/K$.  For $\eta\in K^{\times}$ let $\sigma_{\eta}$
denote the element of $\Gal(K^{ab}/K)$ that corresponds to
$\eta$ under class field theory.  For each $\beta\in K$ let
$\rho_{\beta}\in K^{alg}$ be a root of the polynomial
$X^p-X-\beta$.  Then $\rho_{\beta}\in K^{ab}$, so we can
set
$[\beta,\eta)_K=\sigma_{\eta}(\rho_{\beta})-\rho_{\beta}$.
This construction is independent of the choice of
$\rho_{\beta}$
and defines a bilinear pairing
\begin{equation} \label{pairing}
[\;\:,\;)_K:K\times K^{\times}\lra\Fp.
\end{equation}
Let $\wp:K\ra K$ denote the Artin-Schreier operator
$\wp(x)=x^p-x$.  Then the kernel of $[\;\:,\;)_K$ on the left
is $\wp K$, and the kernel of $[\;\:,\;)_K$ on the right is
$(K^{\times})^p$.  Let $\Omega_K=K\, dt$ denote the module
of K\"ahler differentials of $K$, and for $\omega\in\Omega_K$
let $\Res(\omega)$ denote the residue of $\omega$ at 0.  By
Schmid's theorem \cite[XIV\,\S5, Prop.\,15]{cl} the pairing
(\ref{pairing}) can be computed in terms of the trace
map $\Tr_{\Kbar/\Fp}:\Kbar\ra\Fp$ using the formula
\begin{equation} \label{schmid}
[\beta,\eta)_K=\Tr_{\Kbar/\Fp}\left(\Res
\left(\beta\frac{d\eta}{\eta}\right)\right).
\end{equation}

     Let $B$ denote the subgroup of the additive group
of $K$ which corresponds under Artin-Schreier theory to
the $(\Z/p\Z)^{\nu}$-extension $L/K$:
\begin{equation}
B=\{\beta\in K:X^p-X-\beta\text{ splits over } L\}.
\end{equation}
Recall that $H=\N_{L/K}(L^{\times})$ is the subgroup
of $K^{\times}$ that corresponds to $L/K$ under class field
theory.  The groups $B$ and $H$ are orthogonal
complements of each other with respect to $[\;\:,\;)_K$:
\begin{align}
H&=\left\{\eta\in K^{\times}:[\beta,\eta)_K=0
\text{ for all }\beta\in B\right\} \\
\label{Hcomp}
B&=\left\{\beta\in K:[\beta,\eta)_K=0
\text{ for all }\eta\in H\right\}.
\end{align}

     We need the following lemma regarding the
Artin-Hasse exponential series $E_p(X)\in\Fp[[X]]$:

\begin{lemma} \label{dlogK}
We have $\dst E_p'(X)=E_p(X)\cdot\frac{\lambda(X)}{X}$,
with $\lambda(X)=X+X^p+X^{p^2}+\cdots$.
\end{lemma}

\proof Recall that $E_0(X)\in\Z_p[[X]]$ is the
Artin-Hasse series in characteristic 0.  By
\cite[p.\,52]{pdg} we have
\begin{equation}
E_0(X)
=\exp\left(X+\frac{1}{p}X^p+\frac{1}{p^2}X^{p^2}+\cdots\right),
\end{equation}
where $\exp(X)$ is the usual exponential series.
By taking the formal derivative of this equation we get
\begin{equation} \label{dlogEp}
E_0'(X)=E_0(X)\cdot\frac{\lambda(X)}{X}.
\end{equation}
Since $E_p(X)$ is the image of $E_0(X)$ under
the canonical map $\gamma:\Z_p[[X]]\ra\Fp[[X]]$, the
lemma follows by applying $\gamma$ to (\ref{dlogEp}).
\qed \medskip

     Using Lemma~\ref{dlogK} we see that for
$\alpha\in\M_K$ we have
\begin{equation} \label{dlogE}
\frac{dE_p(\alpha)}{E_p(\alpha)}=
\frac{E_p'(\alpha)\alpha'\,dt}{E_p(\alpha)}=
\frac{\lambda(\alpha)\alpha'\,dt}{\alpha},
\end{equation}
where $\alpha'$ denotes the formal derivative of
$\alpha$ with respect to $t$.  By applying this formula
to (\ref{prod}) we get
\begin{align} \label{dlog}
\frac{d\N_{L/K}(E_p(r\pi_L^k)))}{\N_{L/K}(E_p(r\pi_L^k)))}
\equiv k\lambda(r^{p^{\nu}}t^k)\frac{dt}{t}+
\sum_{j=0}^{\nu-1}\sum_{h\in S_{kj}}
h\lambda(-kc_{kp^j,h}r^{p^j}t^h)\frac{dt}{t}
\!\pmod{\M_K^b\,dt}.
\end{align}

     Let $\Kbar_0$ be a subgroup of the additive group
of $\Kbar$ which is a complement of $\wp\Kbar$; then
$\Kbar_0$ is cyclic of order $p$.  Define
\begin{equation}
K_0=\{x_0+x_1t^{-1}+\dots+x_st^{-s}\in K:s\ge0,\;
x_0\in \Kbar_0,\;x_i\in\Kbar,\;x_{pi}=0\text{ for }i\ge1\}.
\end{equation}
Also define $B_0=B\cap K_0$.  Then $K=K_0\oplus\wp K$
and $B=B_0\oplus\wp K$, so $B/\wp K\cong B_0$.  Let
$\xi=x_0+x_1t^{-1}+\dots+x_bt^{-b}$ be an element of
$K_0\cap\M_K^{-b}$.  We wish to determine the conditions
that $x_0,x_1,\dots,x_b$ must satisfy for $\xi$ to lie
in $B_0$.  Since $B_0\subset K_0\cap\M_K^{-b}$, this
will give a characterization of $B_0$, and hence of $B$.

     Since $B$ is the orthogonal complement of $H$ with
respect to $[\;\:,\;)_K$ we have $\xi\in B_0$ if and
only if $[\xi,\eta)_K=0$ for all $\eta\in H$.  Since
$v_K(\xi)\ge-b$, by (\ref{schmid}) we have
$[\xi,\eta)_K=0$ for all $\eta\in U_K^{b+1}$.  Hence by
Proposition~\ref{generate} we see that
$\xi\in B_0$ if and only if $[\xi,t)_K=0$ and
$[\xi,\N_{L/K}(E_p(r\pi_L^k)))_K=0$ for all $r\in\Kbar$
and all $k$ such that $1\le k\le b$ and $p\nmid k$.
Using (\ref{schmid}) we deduce that $\xi\in B_0$ if and
only if $\Tr_{\Kbar/\Fp}(x_0)=0$ and
\begin{equation} \label{TrRes}
\Tr_{\Kbar/\Fp}\left(\Res\left(\xi\cdot
\frac{d\N_{L/K}(E_p(r\pi_L^k))}{\N_{L/K}(E_p(r\pi_L^k))}
\right)\right)=0
\end{equation}
for all $r\in\Kbar$ and all $1\le k\le b$ such that $p\nmid
k$.  Since $x_0\in \Kbar_0$ and $\wp\Kbar=\ker(\Tr_{\Kbar/\Fp})$,
we have $\Tr_{\Kbar/\Fp}(x_0)=0$ if and only if $x_0=0$.
Using (\ref{dlog}) we see that (\ref{TrRes}) is equivalent to
\begin{equation} \label{Tr}
\Tr_{\Kbar/\Fp}\left(x_kr^{p^{\nu}}-
\sum_{j=0}^{\nu-1}\sum_{h\in S_{kj}}
\,hc_{kp^j,h}x_hr^{p^j}\right)=0.
\end{equation}

     Set $m=[\Kbar:\Fp]$ and let $\tau$ be the Frobenius
automorphism of $\Kbar$.  Then
\begin{equation}
\Gal(\Kbar/\Fp)=\{\tau^0,\tau^1,\dots,\tau^{m-1}\},
\end{equation}
and we can rewrite (\ref{Tr}) in the form
\begin{equation} \label{xkpi}
\sum_{i=0}^{m-1}x_k^{p^i}\tau^{i+\nu}(r)=
\sum_{i=0}^{m-1}\sum_{j=0}^{\nu-1}\sum_{h\in S_{kj}}
hc_{kp^j,h}^{p^i}x_h^{p^i}\tau^{i+j}(r).
\end{equation}
Since (\ref{xkpi}) holds for all $r\in\Kbar$ we get
\begin{equation}
\sum_{i=0}^{m-1}x_k^{p^i}\tau^{i+\nu}=
\sum_{i=0}^{m-1}\sum_{j=0}^{\nu-1}\sum_{h\in S_{kj}}
hc_{kp^j,h}^{p^i}x_h^{p^i}\tau^{i+j}.
\end{equation}
Since $\tau^u=\tau^v$ if and only if $u\equiv v\pmod{m}$
this implies
\begin{equation} \label{xkpij}
\sum_{i=0}^{m-1}x_k^{p^{i-\nu}}\tau^i=
\sum_{i=0}^{m-1}\sum_{j=0}^{\nu-1}\sum_{h\in S_{kj}}
hc_{kp^j,h}^{p^{i-j}}x_h^{p^{i-j}}\tau^i.
\end{equation}
It follows from the $\Kbar$-independence of
$\tau^0,\tau^1,\dots,\tau^{m-1}$ that (\ref{xkpij})
holds if and only if
\begin{equation} \label{xk}
x_k=\sum_{j=0}^{\nu-1}\sum_{h\in S_{kj}}
hc_{kp^j,h}^{p^{\nu-j}}x_h^{p^{\nu-j}}.
\end{equation}

     Let $V$ denote the subgroup of $\Kbar$ such that
$B_0+\M_K^{-b+1}=Vt^{-b}+\M_K^{-b+1}$.  Since $f_{bp^j}=b$
we have $S_{bj}=\{b\}$ for $0\le j\le \nu-1$.  Hence
by (\ref{xk}) every element of $V$ is a root of the
polynomial
\begin{equation} \label{q}
q(X)=X-\sum_{j=0}^{\nu-1}\,bc_{bp^j,b}^{p^{\nu-j}}X^{p^{\nu-j}}.
\end{equation}
Since the only ramification break of $L/K$ is $b$ we
have $V\cong B_0\cong(\Z/p\Z)^{\nu}$.  Hence $V$ is equal to
the set of roots of $q(X)$.  Let $\psi:V\ra B_0$ be
the inverse of the projection of
$B_0$ onto $V$.  Then for each $x_b\in V$ there are uniquely
determined $x_j\in \Kbar$ such that
$\psi(x_b)=x_1t^{-1}+x_2t^{-2}+\dots+x_bt^{-b}$, with $x_i=0$
for $p\mid i$.  For each $k$ the map $\psi_k:V\ra\Kbar$ which
takes $x_b$ to $x_k$ is a group homomorphism.

     Since the maps $\phi_i:V\ra\Kbar$ defined by
$\phi_i(x)=x^{p^i}$ for $1\le i\le\nu$ are linearly
independent over $\Kbar$, they form a basis for the
$\Kbar$-vector space
of homomorphisms from $V$ to $\Kbar$.  It follows that
for $0\le k\le b$ there are uniquely determined
$w_{jk}\in\Kbar$ such that
\begin{equation} \label{psik}
\psi_k(x_b)=w_{0k}x_b^{p^{\nu}}+w_{1k}x_b^{p^{\nu-1}}+\dots
+w_{\nu-1,k}x_b^p
\end{equation}
By setting $w_j=w_{j1}t^{-1}+w_{j2}t^{-2}+\dots+w_{jb}t^{-b}$
we get
\begin{equation} \label{wj}
\psi(x_b)=w_0x_b^{p^{\nu}}+w_1x_b^{p^{\nu-1}}+\dots
+w_{\nu-1}x_b^p
\end{equation}
for all $x_b\in V$.  Hence
\begin{equation} \label{B}
B_0=\{w_0x_b^{p^{\nu}}+w_1x_b^{p^{\nu-1}}+\dots+w_{\nu-1}x_b^p:
x_b\in V\}.
\end{equation}
Of course we have $v_K(w_j)\ge-b$, and $v_K(w_j)\le-1$
if $w_j\not=0$.  In addition, since
$\widehat{\imath}_0=i_0=bp^{\nu}-b$, it follows from
(\ref{bpn}) that $c_{b,b}\not=0$.  Hence $v_K(w_0)=-b$.

     The maps $\psi_k$ can be determined inductively
using (\ref{xk}).  In fact the computations in this
section can be
used to derive an algorithm for computing the indices of
inseparability of the extension $L/K$ in terms of the
norm group $H=\N_{L/K}(L^{\times})$: \medskip

\begin{algorithm} \label{compute}
Given $H$, use (\ref{schmid})
and (\ref{Hcomp}) to determine $B$ and $B_0=B\cap K_0$.
For each $k$ such that $1\le k\le b$ we get the map
$\psi_k:V\ra\Kbar$, from which we obtain
$w_{jk}\in\Kbar$ satisfying (\ref{psik}).  Use
reverse induction on $k$ to determine $c_{kp^j,b}$ for all
$k,j$ such that $0\le k\le b$, $p\nmid k$, and
$1\le j\le\nu-1$: For the base case note that since
$V$ is equal to the set of roots of
$q(X)$, equation (\ref{q}) determines $c_{bp^j,b}$ for all
$j$.  Now let $1\le k<b$ be such that $p\nmid k$ and assume
that we have computed $c_{rp^j,b}$ for all $j,r$ such that
$k<r\le b$ and $p\nmid r$.  If $h\in S_{kj}$ and $h<b$ then
$kp^j<kp^j+(b-h)p^{\nu}<bp^j$, so
$c_{kp^j,h}=c_{kp^j+(b-h)p^{\nu},b}$ has been determined.
Hence (\ref{xk}) can be used to compute $c_{kp^j,b}$ for
$0\le j\le\nu-1$.  Once all the $c_{kp^j,b}$ have been
found use (\ref{bpn}) to compute the indices of
inseparability of $L/K$;
since $i_j\ge bp^{\nu}-bp^j$ this computation does not
depend on $c_{kp^j,b}$ for $k>b$.
\end{algorithm}

\section{Some explicit formulas} \label{explicit}

In this section we continue to assume that $K$ is a local
field of characteristic $p$ with finite residue field
$\Kbar$, and that $L/K$ is a totally ramified elementary
abelian $p$-extension
of degree $p^{\nu}$ with a single ramification break $b\ge1$.
We use the results of sections \ref{Ab} and \ref{AS} to prove
theorems which relate the indices of inseparability of $L/K$
to the description of $L/K$ in terms of Artin-Schreier
theory.  Our approach is to use (\ref{xk}), (\ref{B}), and
(\ref{bpn}) to get relations between the $K$-valuations of
$w_0,w_1,\dots,w_{\nu-1}$ and the indices of inseparability
$i_0,i_1,\dots,i_{\nu}$.  In some
cases we are able to derive explicit formulas for $i_j$ in
terms of the valuations of the $w_i$.  We then use
Schmid's formula (\ref{schmid})
to translate theorems expressed in terms of
Artin-Schreier theory into theorems expressed in terms
of local class field theory.
We remark that although $B_0$ and $w_0,w_1,\dots,w_{\nu-1}$
depend on $t$, and hence on the choice of $\pi_L$, the results
of this section do not, of course, depend on this choice.

\begin{theorem} \label{bp}
Let $L/K$ be a totally ramified $(\Z/p\Z)^{\nu}$-extension
with a single ramification break $b$ satisfying
$1\le b\le p-1$.  Then for $0\le j\le\nu-1$ we have
\begin{equation} \label{ijformula}
i_j=bp^{\nu}+\min\{p^{j'}v_K(w_{j'}):0\le j'\le j\}.
\end{equation}
\end{theorem}

\proof Since $b<p$ we have $S_{kj}=\{b\}$ for
$1\le k\le b$ and $0\le j\le\nu-1$.  Hence for
$1\le k\le b$ formula (\ref{xk}) simplifies to
\begin{equation} \label{simp}
\psi_k(x_b)
=\sum_{j=0}^{\nu-1}\,bc_{kp^j,b}^{p^{\nu-j}}x_b^{p^{\nu-j}}.
\end{equation}
Comparing this with (\ref{psik}) and (\ref{wj})
we get
\begin{equation} \label{simpwj}
w_j=\sum_{k=1}^b\,bc_{kp^j,b}^{p^{\nu-j}}t^{-k}
\end{equation}
for $0\le j\le\nu-1$.  

     We now prove our claim by induction on $j$.
Since $i_0=bp^{\nu}-b$ and $v_K(w_0)=-b$ the claim holds
for $j=0$.  Let $1\le j\le\nu-1$ and assume that the
claim holds for $j-1$.
Suppose $w_j=0$.  Then $c_{kp^j,b}=0$ for $1\le k\le b$.
Since $b<p$ and $i_{j-1}\ge bp^{\nu}-bp^{j-1}$, by
(\ref{bpn}) we get $i_j=i_{j-1}$.  We also have
\begin{equation}
\min\{p^{j'}v_K(w_{j'}):0\le j'\le j\}=
\min\{p^{j'}v_K(w_{j'}):0\le j'\le j-1\}.
\end{equation}
Therefore the claim for $j$ follows from the claim for
$j-1$ in this case.  If $w_j\not=0$ then $v_K(w_j)=-k$
for some $k$ such that $1\le k\le b$.  Since $b<p$ this
implies that the right side of (\ref{ijformula}) is
equal to $bp^{\nu}-kp^j$.  By (\ref{simpwj}) we have
$c_{kp^j,b}\not=0$ and $c_{lp^j,b}=0$ for all $l$ such
that $k<l\le b$.  Since $i_{j-1}\ge bp^{\nu}-bp^{j-1}$
it follows from (\ref{bpn}) that $i_j=bp^{\nu}-kp^j$.
Hence the claim holds for $j$. \qed \medskip

     Although there does not seem to be a simple formula
for the indices of inseparability similar to
Theorem~\ref{bp} which is valid for all $b$, it is
possible to get some information about $i_j$ in the
general case.  Let $\Fpn$ denote the unique subfield of
$\Kbar^{alg}$ with $p^{\nu}$ elements.

\begin{lemma} \label{B0}
Let $L/K$ be a totally ramified
$(\Z/p\Z)^{\nu}$-extension with a single ramification
break $b\ge1$, and let
$0\le k\le b-1$.  Then the following are equivalent:
\begin{enumerate}
\item $\Fpn\subset\Kbar$ and $B_0+\M_K^{-k}$ is an
$\Fpn$-subspace of $\M_K^{-b}$.
\item $w_j\in\M_K^{-k}$ for $1\le j\le\nu-1$.
\end{enumerate}
\end{lemma}

\proof Suppose condition 2 holds for $k$.  Then
$w_{jb}=0$ for $1\le j\le\nu-1$, so we have
$x_b=\psi_b(x_b)=w_{0b}x_b^{p^{\nu}}$.  Since
$V\subset\Kbar$ is the set of solutions to this equation
we deduce that $\Fpn\subset\Kbar$ and $V$ is a vector
space over $\Fpn$.  Furthermore, by (\ref{B}) we have
\begin{equation}
B_0+\M_K^{-k}=\{w_0x_b^{p^{\nu}}:x_b\in V\}+\M_K^{-k},
\end{equation}
so $B_0+\M_K^{-k}$ is an $\Fpn$-subspace of $\M_K^{-b}$.
Conversely, if condition 1 holds for $k$ then $V$ is a
vector space over $\Fpn$, and for every $c\in\Fpn$ and
$x_b\in V$ we have
\begin{equation}
\psi(cx_b)\equiv c\psi(x_b)\pmod{\M_K^{-b+1}}.  
\end{equation}
Since $(B_0+\M_K^{-k})/\M_K^{-k}$ is a one-dimensional
vector space over $\Fpn$ this implies
\begin{equation}
\psi(cx_b)\equiv c\psi(x_b)\pmod{\M_K^{-k}}.  
\end{equation}
It follows that $\psi_i(cx_b)=c\psi_i(x_b)$ for
$k<i\le b$, so we have $w_{ji}=0$ for $1\le j\le\nu-1$
and $k<i\le b$.  Therefore condition 2 holds for $k$.
\qed

\begin{theorem} \label{Zpn}
Let $L/K$ be a totally ramified
$(\Z/p\Z)^{\nu}$-extension with a single ramification
break $b\ge1$, and let $k$ be an integer such that
$\lceil b/p\rceil\le k\le b-1$.  Then the following are
equivalent:
\begin{enumerate}
\item $i_j\ge bp^{\nu}-kp^j$ for $1\le j\le\nu-1$.
\item $\Fpn\subset\Kbar$ and $B_0+\M_K^{-k}$ is an
$\Fpn$-subspace of $\M_K^{-b}$.
\item $w_j\in\M_K^{-k}$ for $1\le j\le\nu-1$.
\end{enumerate}
\end{theorem}

\proof The equivalence of conditions 2 and 3 follows
from Lemma~\ref{B0}.  To prove the equivalence of
conditions 1 and 3 we
use reverse induction on $k$.  By (\ref{bpn}) we see
that condition 1 holds for $k=b-1$ if and only if
$c_{bp^j,b}=0$ for $1\le j\le\nu-1$.  Using (\ref{q}) we
see that this is equivalent to
$q(X)=X-c_{b,b}^{p^{\nu}}X^{p^{\nu}}$, which holds if
and only if condition 3 holds for $k=b-1$.  This proves
the base case.  Now let
$\lceil b/p\rceil<r\le b-1$ and assume that
the theorem holds for $k=r$.  If $p\mid r$ then
condition 3 holds for $k=r$ if and only if condition 3
holds for $k=r-1$, and condition 1 holds for $k=r$ if
and only if condition 1 holds for $k=r-1$ (since
$p^{j+1}\nmid i_j$).  Hence the theorem holds for $k=r-1$
in this case.  From now on we assume that the theorem
holds for $k=r$ with $\lceil b/p\rceil<r\le b-1$ and
$p\nmid r$.

     Suppose condition 1 holds for $k=r-1$.  Then by
(\ref{without}) we see that $c_{rp^j,h}=0$ for all
$(j,h)$ such that $1\le j\le\nu-1$ and $1\le h\le b$.
Hence by (\ref{xk}) we have
\begin{equation} \label{xr}
\psi_r(x_b)=\sum_{h\in S_{r0}}\,
hc_{r,h}^{p^{\nu}}\psi_h(x_b)^{p^{\nu}}.
\end{equation}
Since condition 1 holds for $k=r-1$ it also holds for
$k=r$.  It follows from the inductive assumption that
condition 3 holds for $k=r$.  Hence for every $h$ such
that $r<h\le b$ we have
$\psi_h(x_b)=w_{0h}x_b^{p^{\nu}}$ for all $x_b\in V$.  In
particular, $x_b=w_{0b}x_b^{p^{\nu}}$ with
$w_{0b}\not=0$.  It follows that $\psi_h(x_b)=w_{0h}w_{0b}^{-1}x_b$
for $r<h\le b$.  Since $r<b$, every $h\in S_{r0}$
satisfies $r<h\le b$.  Hence by (\ref{xr}) we have
\begin{equation} \label{xk1}
\psi_r(x_b)=\sum_{h\in S_{r0}}\,
hc_{r,h}^{p^{\nu}}w_{0h}^{p^{\nu}}w_{0b}^{-p^{\nu}}x_b^{p^{\nu}}.
\end{equation}
Since $v_K(w_j)\ge-r$ for $1\le j\le\nu-1$, and
(\ref{xk1}) expresses $\psi_r(x_b)$ as a
$\Kbar$-multiple of $x_b^{p^{\nu}}$, we deduce that
$v_K(w_j)\ge-r+1$ for $1\le j\le\nu-1$.  Therefore
condition 3 holds for $k=r-1$.

     Assume conversely that condition 3 holds for
$k=r-1$.  Then condition 3 also holds for $k=r$, so by
the inductive assumption, condition 1 holds for $k=r$.
It follows from (\ref{without}) that $c_{rp^j,h}=0$ for
all $(j,h)$ such that $1\le j\le\nu-1$ and $1\le h<b$.
Therefore by (\ref{xk}) we get
\begin{equation}
\psi_r(x_b)=\sum_{j=1}^{\nu-1}\,
bc_{rp^j,b}^{p^{\nu-j}}x_b^{p^{\nu-j}}
+\sum_{h\in S_{r0}}\,hc_{r,h}^{p^{\nu}}\psi_h(x_b)^{p^{\nu}}.
\end{equation}
As in the preceding paragraph we have
$\psi_h(x_b)=w_{0h}w_{0b}^{-1}x_b$ for every $h\in S_{r0}$.
Hence
\begin{equation} \label{xrsimp}
\psi_r(x_b)=\sum_{j=1}^{\nu-1}\,
bc_{rp^j,b}^{p^{\nu-j}}x_b^{p^{\nu-j}}+\sum_{h\in S_{r0}}\,
hc_{r,h}^{p^{\nu}}w_{0h}^{p^{\nu}}w_{0b}^{-p^{\nu}}x_b^{p^{\nu}}.
\end{equation}
It follows from (\ref{psik}) that
$w_{jr}=bc_{rp^j,b}^{p^{\nu-j}}$ for $1\le j\le\nu-1$.
Hence by condition 3 for $k=r-1$ we get $c_{rp^j,b}=0$
for $1\le j\le\nu-1$.  Suppose condition 1 does not hold
for $k=r-1$, and let $1\le j\le\nu-1$ be minimum such
that $i_j<bp^{\nu}-(r-1)p^j$.  Since $r-1\ge b/p$ we
have $i_1<i_0$, and for $2\le j\le\nu-1$ we have
$i_j<i_{j-1}$ by the minimality of $j$.  Hence
$p^j\mid i_j$ for $1\le j\le\nu-1$.
Since condition 1 holds for $r$ it follows that
$i_j=bp^{\nu}-rp^j$.  By (\ref{bpn}) this implies
$c_{rp^j,b}\not=0$, a contradiction.  Hence condition 1
holds for $k=r-1$.  \qed

\begin{cor}
Let $L/K$ be a totally ramified
$(\Z/p\Z)^{\nu}$-extension with a single ramification
break $b$ such that $b\ge2$ and $\Fpn$ is not
contained in $\Kbar$.  Then there is at least one
$1\le j\le\nu-1$ such that $i_j=bp^{\nu}-bp^j$.
\end{cor}

\proof Since $\Fpn$ is not contained in $\Kbar$,
condition 2 in Theorem~\ref{Zpn} does not hold for
$k=b-1$.  Since $b\ge2$ we have $b-1\ge\lceil
b/p\rceil$, so condition 1 also fails for $k=b-1$.
Let $1\le j\le\nu-1$ be minimum such that
$i_j<bp^{\nu}-(b-1)p^j$.  As in the proof of
Theorem~\ref{Zpn} we get $p^j\mid i_j$.
Since we also have $i_j\ge bp^{\nu}-bp^j$ we conclude
that $i_j=bp^{\nu}-bp^j$. \qed

\begin{cor} \label{Zp2}
Let $L/K$ be a totally ramified $(\Z/p\Z)^2$-extension
with a single ramification break $b\ge1$.  Then
\begin{align}
i_1&=bp^2+\min\{v_K(w_0),pv_K(w_1)\} \\
&=bp^2+\min\{-b,pv_K(w_1)\}.
\end{align}
\end{cor}

\proof Set $l=-v_K(w_1)$ and $d=\lceil b/p\rceil$.  If
$d+1\le l\le b$ then by Theorem~\ref{Zpn} we have
\begin{equation}
bp^2-lp\le i_1<bp^2-(l-1)p<bp^2-b=i_0.
\end{equation}
Since $i_1<i_0$ we have $p\mid i_1$, and hence
$i_1=bp^2-lp$.  If $l\le d$ then the conditions of
Theorem~\ref{Zpn} hold for $k=d$, so we have
$i_1\ge bp^2-dp$.  Furthermore, (\ref{xrsimp}) is valid
with $\nu=2$ and $r=d$, so $w_{1d}=bc_{dp,b}^p$.
If $l=d$ we get $c_{dp,b}\not=0$, and hence
$i_1=bp^2-dp$ by (\ref{bpn}), while if $l<d$
then $c_{dp,b}=0$, and hence $i_1=bp^2-b$.  We conclude
that $i_1=bp^2+\min\{-b,pv_K(w_1)\}$ in every case.
\qed \medskip

     The power series $E_p(X)$ induces a bijection
from $\M_K$ onto $U_K^1$.  Let $\Lambda_K:U_K^1\ra\M_K$
denote the inverse of this bijection.  Recall that
$H=\N_{L/K}(L^{\times})$ is the subgroup of $K^{\times}$
which corresponds to $L/K$ under class field theory.

\begin{lemma} \label{subgroup}
Suppose $i\ge\lceil b/p\rceil$.  Then
$\Lambda_K(H\cap U_K^i)$ is a subgroup of the additive
group of $\M_K$.
\end{lemma}

\proof By equation (6) in \cite[p.\,52]{pdg}, for
$\alpha_1,\alpha_2\in\M_K^i$ we have
\begin{equation} \label{Epcong}
E_p(\alpha_1+\alpha_2)\equiv
E_p(\alpha_1)E_p(\alpha_2)\pmod{\M_K^{pi}}.
\end{equation}
Since $E_p(\M_K^i)=U_K^i$, it follows that
\begin{equation} \label{Lambdamult}
\Lambda_K(u_1u_2)\equiv
\Lambda_K(u_1)+\Lambda_K(u_2)\pmod{\M_K^{pi}}
\end{equation}
for $u_1,u_2\in U_K^i$.  Since $pi\ge b+1$ we have
$U_K^{pi}\subset U_K^{b+1}\subset H$, and hence
$\M_K^{pi}\subset\Lambda_K(H\cap U_K^i)$.  Thus
$\Lambda_K(H\cap U_K^i)$ is a subgroup of $\M_K$. \qed

\begin{lemma} \label{linear}
Let $\alpha\in\M_K^i$ with $i\ge\lceil b/p\rceil$.  Then
for every $\beta\in\M_K^{-b}$ and $\zeta\in\Kbar$ we have
\begin{equation}
[\zeta\beta,E_p(\alpha))_K
=[\beta,E_p(\zeta\alpha))_K.
\end{equation}
\end{lemma}

\proof Since $pi-1\ge b$ it follows from (\ref{dlogE})
that
\begin{alignat}{2}
\frac{dE_p(\alpha)}{E_p(\alpha)}&\equiv
\alpha'\,dt&&\pmod{\M_K^b\, dt} \\
\frac{dE_p(\zeta\alpha)}{E_p(\zeta\alpha)}
&\equiv\zeta\alpha'\,dt&&\pmod{\M_K^b\, dt}.
\end{alignat}
The lemma now follows from Schmid's formula
(\ref{schmid}). \qed

\begin{theorem} \label{sub}
Let $k\ge\lceil b/p\rceil-1$.  Then $\Lambda_K(H\cap
U_K^{k+1})$ is an $\Fpn$-subspace of $\M_K$ if and only if
$B_0+\M_K^{-k}$ is an $\Fpn$-subspace of $K$.
\end{theorem}

\proof If $\Fpn$ is not contained in $\Kbar$ then
neither $\Lambda_K(H\cap U_K^{k+1})$ nor $B_0+\M_K^{-k}$
is a vector space over $\Fpn$.  If $\Fpn\subset\Kbar$
and $k\ge b$ then $\Lambda_K(H\cap
U_K^{k+1})=\M_K^{k+1}$ and $B_0+\M_K^{-k}=\M_K^{-k}$
are both vector spaces over $\Fpn$.  Hence we may
assume $\Fpn\subset\Kbar$ and $k\le b-1$.

     Suppose $B_0+\M_K^{-k}$ is an $\Fpn$-subspace of
$K$.  By Lemma~\ref{subgroup} we
see that $\Lambda_K(H\cap U_K^{k+1})$ is a
subgroup of $\M_K$.  We need to show
that $\Lambda_K(H\cap U_K^{k+1})$ is stable under
multiplication by elements of $\Fpn$.
Let $\beta\in B_0+\M_K^{-k}$ and
$\eta\in H\cap U_K^{k+1}$.  Then $\eta=E_p(\alpha)$
for a uniquely determined $\alpha\in\M_K^{k+1}$.  For
$\zeta\in\Fpn$ we have $\zeta\beta\in B_0+\M_K^{-k}$.
Since
\begin{equation}
[B_0,H)_K=[\M_K^{-k},U_K^{k+1})_K=0,
\end{equation}
by Lemma~\ref{linear} we get
$[\beta,E_p(\zeta\alpha))_K=0$.  Hence
$E_p(\zeta\alpha)$ lies in the orthogonal complement
of $B_0+\M_K^{-k}$ with respect to $[\;\:,\;)_K$.  Since
$B=B_0\oplus\wp K$, the orthogonal complement of $B_0$ is
equal to the orthogonal complement of $B$, which is $H$.
Since we also have $E_p(\zeta\alpha)\in U_K^{k+1}$
we get $E_p(\zeta\alpha)\in H\cap U_K^{k+1}$.  Hence
$\zeta\alpha\in\Lambda_K(H\cap U_K^{k+1})$ for all
$\zeta\in\Kbar$, so $\Lambda_K(H\cap U_K^{k+1})$ is an
$\Fpn$-subspace of $\M_K$.

     Conversely, suppose that $\Lambda_K(H\cap U_K^{k+1})$
is an $\Fpn$-subspace of $\M_K$.  Then
for every $\eta=E_p(\alpha)$ in $H\cap U_K^{k+1}$ and
every $\zeta\in\Fpn$ we have $E_p(\zeta\alpha)\in H\cap
U_K^{k+1}$.  Therefore by Lemma~\ref{linear} we have
\begin{equation}
[\zeta\beta,E_p(\alpha))_K
=[\beta,E_p(\zeta\alpha))_K=0
\end{equation}
for every $\beta\in B_0$.
Hence $\zeta\beta$ lies in the orthogonal complement of
$H\cap U_K^{k+1}$.  Since $H$ and
$U_K^{k+1}(K^{\times})^p$ are closed finite-index
subgroups of $K^{\times}$ we have
\begin{align}
(H\cap U_K^{k+1})^{\perp}
&=((H\cap U_K^{k+1})(K^{\times})^p)^{\perp} \\
&=(H\cap(U_K^{k+1}(K^{\times})^p))^{\perp} \\
&=H^{\perp}+(U_K^{k+1}(K^{\times})^p)^{\perp} \\
&=B_0+\M_K^{-k}+\wp K.
\end{align}
Since we also have
$\zeta\beta\in(K_0\cap\M_K^{-b})+\Kbar$ we get
\begin{equation}
\zeta\beta\in(B_0+\M_K^{-k}+\wp K)\cap
((K_0\cap\M_K^{-b})+\Kbar)\subset B_0+\M_K^{-k}.
\end{equation}
It follows that $B_0+\M_K^{-k}$ is an $\Fpn$-subspace of
$K$. \qed \medskip

     By combining Theorem~\ref{sub} with
Theorem~\ref{Zpn} we get the following result:

\begin{cor} \label{Lp}
Let $L/K$ be a totally ramified
$(\Z/p\Z)^{\nu}$-extension with a single ramification
break $b\ge1$ and let $k$ be an integer such that
$\lceil b/p\rceil\le k\le b-1$.  Then the following are
equivalent:
\begin{enumerate}
\item $i_j\ge bp^{\nu}-kp^j$ for $1\le j\le\nu-1$.
\item $\Fpn\subset\Kbar$ and $\Lambda_K(H\cap
U_K^{k+1})$ is an $\Fpn$-subspace of $\M_K$.
\end{enumerate}
\end{cor}

     We also have the following reformulation of
Corollary~\ref{Zp2}:

\begin{cor} \label{combining}
Let $L/K$ be a totally ramified $(\Z/p\Z)^2$-extension
with a single ramification break $b\ge1$.
If $\Fps\not\subset\Kbar$ let
$k=b$; otherwise let $k$ be the smallest nonnegative
integer such that $\Lambda_K(H\cap U_K^{k+1})$ is an
$\Fps$-subspace of $\M_K$.  Then
$i_1=p^2b-\max\{b,pk\}$.
\end{cor}

\proof If $\Fps\not\subset\Kbar$ then $B_0+\M_K^{-b+1}$
is not a vector space over $\Fps$, so $v_K(w_1)=-b$
by Lemma~\ref{B0}.  Hence by Corollary~\ref{Zp2} we
have $i_1=bp^2-bp$.  If $\Fps\subset\Kbar$ set
$d=\lceil b/p\rceil$.  If $k\ge d$ then by
Theorem~\ref{sub} we see that $k$ is the smallest
nonnegative integer such that $B_0+\M_K^{-k}$ is an
$\Fps$-subspace of $\M_K^{-b}$.  Hence by Lemma~\ref{B0} 
we have $v_K(w_1)=-k$, so by Corollary~\ref{Zp2}
we get $i_1=p^2b-pk$.  If $k<d$ then
$\Lambda_K(H\cap U_K^d)$ is an $\Fps$-subspace of
$\M_K$.  Hence by Theorem~\ref{sub} we see that
$B_0+\M_K^{-d+1}$ is an $\Fps$-subspace of $\M_K^{-b}$.
Hence $v_K(w_1)\ge-d+1$, so by Corollary~\ref{Zp2} we
have $i_1=p^2b-b$.  We conclude that
$i_1=p^2b-\max\{b,pk\}$ in every case. \qed

\section{The case $\char(K)=0$} \label{char0}

In this section we show how to extend the results of
the previous section to apply to certain extensions of
local fields of characteristic 0.  Let $K$ be a finite
extension of the $p$-adic field $\Q_p$ and set
$e=v_K(p)$.  Let $L/K$ be a totally ramified Galois
extension with a single ramification break $b$ such that
$1\le b<e$.  Then $\Gal(L/K)\cong(\Z/p\Z)^{\nu}$ for
some $\nu\ge1$.

     Let $\pi_K$, $\pi_L$ be uniformizers for $K$, $L$
and recall that $R$ denotes the set of Teichm\"uller
representatives for $K$.  Choose $c_h\in R$ such that
$\dst\pi_K=\sum_{h=0}^{\infty}\,c_h\pi_L^{h+p^{\nu}}$
and let $\cbar_h$ denote the image of $c_h$ in $\Kbar$.
Since $K'=\Kbar((T))$ is a local field of characteristic
$p$ with residue field $\Kbar$, by \cite[Th.\,2.2]{heier}
there is a totally ramified extension $L'/K'$ and a
uniformizer $\pi_{L'}$ for $L'$ such that
$\dst T=\sum_{h=0}^{\infty}\,\cbar_h\pi_{L'}^{h+p^{\nu}}$.

     For $0\le j\le\nu-1$ the inequalities
$i_{j+1}\ge bp^{\nu}-bp^{j+1}$ and $i_j\le
i_0=bp^{\nu}-b$ imply that $i_j-i_{j+1}\le
bp^{j+1}-b<ep^{\nu}$.  Hence $i_j=\tilde{\imath}_j$ for
$0\le j\le\nu$.  Since $c_h=0$ if and only if
$\cbar_h=0$, this implies that the extension $L'/K'$ has
the same indices of inseparability as $L/K$.  It follows
that $L'/K'$ has the same ramification data as $L/K$, so
$L'/K'$ has a single ramification break $b$.  It
can be shown using \cite[Prop.\,6.3b]{heier} that
$L'/K'$ is Galois.  Hence
$\Gal(L'/K')\cong(\Z/p\Z)^{\nu}$, so
the results of section~\ref{explicit} apply to
$L'/K'$.  Since $b+1\le e$ there is an isomorphism
\begin{equation}
\rho:\OO_L/\M_L^{(b+1)p^{\nu}}\lra
\OO_{L'}/\M_{L'}^{(b+1)p^{\nu}}
\end{equation}
of $\Kbar$-algebras such that
$\rho(\pi_L+\M_L^{(b+1)p^{\nu}})
=\pi_{L'}+\M_{L'}^{(b+1)p^{\nu}}$.  It follows from the
series expressions for $\pi_K$ and $T$ that
$\rho(\pi_K+\M_L^{(b+1)p^{\nu}})=T+\M_{L'}^{(b+1)p^{\nu}}$.
Set $H=\N_{L/K}(L^{\times})$ and
$H'=\N_{L'/K'}((L')^{\times})$.  Let $\tilde{H}$ denote the
image of $\N_{L/K}(\OO_L^{\times})$ in
\begin{equation}
(\OO_K/\M_K^{b+1})^{\times}\cong
\OO_K^{\times}/U_K^{b+1}
\end{equation}
and let $\tilde{H}'$ denote the
image of $\N_{L'/K'}(\OO_{L'}^{\times})$ in
\begin{equation}
(\OO_{K'}/\M_{K'}^{b+1})^{\times}\cong
\OO_{K'}^{\times}/U_{K'}^{b+1}.
\end{equation}
Since the norm of
$u\in\OO_L^{\times}$ is the determinant of the
$\OO_K$-linear map from $\OO_L$ to $\OO_L$ defined by
$x\mapsto ux$, we have $\rho(\tilde{H})=\tilde{H}'$.

     As before we let $E_0(X)\in\Z_p[[X]]$ denote the
Artin-Hasse exponential series in characteristic 0,
and let $\Lambda_K:U_K^1\ra\M_K$ denote the inverse of
the bijection from $\M_K$ to $U_K^1$ induced by
$E_0(X)$.  Similarly, we have $E_p(X)\in\Fp[[X]]$ and
$\Lambda_{K'}:U_{K'}^1\ra\M_{K'}$.
Let $u\in U_K^1$ and $u'\in U_{K'}^1$ be such that
\begin{equation}
\rho(u+\M_L^{(b+1)p^{\nu}})
=u'+\M_{L'}^{(b+1)p^{\nu}}.
\end{equation}
Since $E_p(X)$ is the image of $E_0(X)\in\Z_p[[X]]$ in
$\Fp[[X]]$ we have
\begin{equation}
\rho(\Lambda_K(u)+\M_L^{(b+1)p^{\nu}})
=\Lambda_{K'}(u')+\M_{K'}^{(b+1)p^{\nu}}.
\end{equation}
Since $\rho(\tilde{H})=\tilde{H}'$, it follows that for
$1\le k\le b$, $\rho$ induces a $\Kbar$-isomorphism from
$\Lambda_K(H\cap U_K^{k+1})/\M_K^{b+1}$ onto
$\Lambda_{K'}(H'\cap U_{K'}^{k+1})/\M_{K'}^{b+1}$.
Therefore $\Lambda_K(H\cap U_K^{k+1})/\M_K^{b+1}$ is an
$\Fpn$-subspace of $\M_K/\M_K^{b+1}$ if and only if
$\Lambda_{K'}(H'\cap U_{K'}^{k+1})/\M_{K'}^{b+1}$ is
an $\Fpn$-subspace of $\M_{K'}/\M_{K'}^{b+1}$.

     Assume $\Fpn\subset\Kbar$ and let
$\Zpn$ denote the ring of integers of the unique
unramified subextension of degree $\nu$ of $K/\Q_p$.
Then $\Lambda_K(H\cap U_K^{k+1})/\M_K^{b+1}$ is an
$\Fpn$-subspace of $\M_K/\M_K^{b+1}$ if and only if
$\Lambda_K(H\cap U_K^{k+1})$ is a $\Zpn$-submodule of
$\M_K$.  Hence we have the following analogs of
Corollary~\ref{Lp} and Corollary~\ref{combining}:

\begin{cor} \label{zero}
Let $K$ be a finite extension of $\Q_p$ and let
$e=v_K(p)$ be the absolute ramification index of $K$.
Let $L/K$ be a totally ramified $(\Z/p\Z)^{\nu}$-extension
with a single ramification break $b$ such that $1\le b<e$,
and set $H=\N_{L/K}(L^{\times})$.  Let $k$ be an integer
such that $\lceil b/p\rceil\le k\le b-1$.  Then the
following are equivalent:
\begin{enumerate}
\item $i_j\ge bp^{\nu}-kp^j$ for $1\le j\le\nu-1$.
\item $\Fpn\subset\Kbar$ and
$\Lambda_K(H\cap U_K^{k+1})$ is a $\Zpn$-submodule
of $\M_K$.
\end{enumerate}
\end{cor}

\begin{cor}
Let $K$, $L$, $e$, $b$ be as in Corollary~\ref{zero} and
assume that $\nu=2$.  If $\Kbar$ does not have a
subfield isomorphic to $\Fps$ let $k=b$; otherwise let
$k$ be the smallest nonnegative integer such that
$\Lambda_K(H\cap U_K^{k+1})$ is a $\Zps$-submodule of
$\M_K$.  Then $i_1=p^2b-\max\{b,pk\}$.
\end{cor}

\end{document}